\documentclass[10pt,a4paper]{article}
\usepackage{color,amsmath,amssymb, lscape,graphicx,ifthen}
\usepackage[col
orlinks=true,
linkcolor=webgreen,
filecolor=webbrown,
citecolor=webgreen]{hyperref}
\definecolor{webgreen}{rgb}{0,.5,0}\definecolor{webbrown}{rgb}{.6,0,0}

\font\meinfont=cmsl8 scaled \magstep1
\def\eg{{\it e.g.},\,}

\def\ie{{\it i.e.},\,}

\def\sspp{\,+\,}
\def\sspm{\,-\,}
\def\sspeq{\,=\,}
\def\sspdef{\, :=\,}
\def\ssptimes{\,\times\,}

\def\sspgeq{\,\geq\, }

\def\sspequiv{\,\equiv\,}
\def\sspin{\,\in\,}

\def\pn{\par\noindent}

\def\pbn{\par\bigskip\noindent}
\def\psn{\par\smallskip\noindent}
\def\Beq{\begin{equation}}
\def\Eeq{\end{equation}}
\def\Beqarray{\begin{eqnarray}}
\def\Eeqarray{\end{eqnarray}}

\normalbaselines
\begin{document}
\bibliographystyle{unsrt}
\rightline{October 31 2023}
\vbox {\vspace{6mm}}
\begin{center}
{\Large {\bf Four Sequences of Length 28 and the Gregorian Calendar}}\\ [9mm]
Wolfdieter L a n g \footnote{ 
\href{mailto:wolfdieter.lang@partner.kit.edu}{\tt wolfdieter.lang@partner.kit.edu},\quad 
\url{http://www.itp.kit.edu/~wl}
                                          } \\[3mm]
\end{center}
\begin{abstract}

It is shown that each sequence giving the number of times a given day of the month falls on a certain day of the week for $400$ successive years of the Gregorian cycle can be composed of various pieces of various length of one of $4$ sequences of length $28$, used periodically.
\end{abstract}
\section{Introduction}
\hskip 1cm
The formula which determines the day of the week, $D \sspeq 0,1,...,6$ for Sunday, Monday, ..., Saturday, for an admissible $d$th day of the month $d\sspin \{1,\,2,\,...,\,31\}$ of a month $m\sspin\in \{1,\,2,\,...,\,12\}$ in a year $y$ of the {\sl Gregorian} calendar \cite{Clavius}, \cite{Wiki1} (in use since October $15$, $1582$, a Friday) is well known. \eg\,\cite{NZ}, II, pp. 357 - 358. 
\begin{eqnarray}
D(d,m,y) &\sspeq &\Bigl{(}d\sspp \lfloor 2.6\, M(m)\sspm 0.2  \rfloor \sspp T(y) \sspp\left  \lfloor {\frac{T(y)}{4}}\right\rfloor\sspp \left \lfloor {\frac{H(y)}{4}}\right \rfloor \nonumber \\
&  & \sspm 2\,H(y)\sspm (L(y)\sspp 1)\,\left \lfloor {\frac{M(m)}{11}}\right \rfloor\Bigr{)}\,{\rm mod}\,7\,. \label{(1.1)}                                                         \end{eqnarray}     
Here $\lfloor . \rfloor$ is the floor function, $M(m)\sspeq (m+10)\,$mod\,$12$ if $M(m)\neq 0$ and $M(2)\sspeq 12$.                                  
$T(y)$ are the last two digits of the year $y$ (the 'tenth') and $H(y)$ are the first two digits of $y\in\{1583,\,...,\,9999\}$. $L(y)\sspeq 1$ if $y$ is a leap year and $L(y)\sspeq 0$ otherwise. In the {\sl Gregorian} calendar $L(y)\sspeq 0$ if and only if $y\,$ mod\,$ 4\,\neq \,0$ or $y\,$mod\,$ 100\sspeq 0$ but  $y\,$mod\,$ 400\,\neq\, 0$. All other years are leap years. {\it E.g.}  $L(1900)\sspeq 0$ but $L(2000)\sspeq1$.
\pn
For the following application it is necessary to exclude nonsense days of a month which are: $d\sspeq 29$ if $m\sspeq 2$ and $L(y)\sspeq 0$, $d\sspeq30$ if $m\sspeq 2$, and $d\sspeq 31$ if $m\sspeq 2,\,4,\,6,\,9,\,11$. Formula ~\ref{(1.1)} satisfies $D(d,m,y)\sspeq D(d,m,y+400)$. Therefore, it suffices to consider $y\, mod\, 400$, {\it i.e.} $y\in [0,399]$ (keeping the start of the {\it Gregorian} calendar in mind). In the formula $y\sspeq 0$ (for $0000$)uses $ H(y)\sspeq 0$ and $T(y)\sspeq 0$.                            
\pn                                                
Every {\sl Gregorian} cycle of $400$ consecutive years has $146\,097 = 7\cdot                      
20871$ days. Such a cycle comprises $97$ leap years and $303$ ordinary years.                      
Any day of the month $d\sspin \{1,\,2,\,...,\,28\}$ appears $12\cdot 400\sspeq 4800$ times in such a {\sl Gregorian} cycle. Because $4800\,$mod$\,7\, \neq\, 0$ such a day of the month cannot be equally distributed over the days of the week. It is well known \cite{Brown}, \cite{Wiki2} that {\it e.g.} $d\sspeq 13$ (hence $d\sspeq 6$ or $20$ or $27$) happens to fall more often on a Friday than on any other day of the week. For these days of the month the distribution is $[687,685,685,687,684,688,684]$ for the days of the week ordered from left to right for Sunday to Saturday. The mean value is $4800/7\simeq 685.71$. Similarly, the first of a month, $d\sspeq 1$ (or $d\sspeq 8,15,22$), falls predominantly on Sunday. The distribution for these days is $[688,684,687,685,685,687,684]$, which is a cyclic permutation of the above given numbers because $d\sspeq 6$ on a Sunday corresponds to $d\sspeq 1$ on a Tuesday, {\it etc}.) 
\pn
$d\sspeq 29$ appears $12\cdot 97 \sspp 11\cdot 303\sspeq 4497$ times during a $400$ year cycle. For $d\sspeq 30$ this number is $11\cdot 400\sspeq 4400$, and for $d\sspeq 31$ it is $7\cdot 400\sspeq 2800$  Even though $2800$ is divisible by $7$ the distribution over the days of the week for day $31$ is also not uniform. 
\pn See {\it Table 1} for the well known (see \eg \cite{Wiki3}) distribution of the days of the month on the days of a week. Each day $d\sspin \{8,\,9,\,...,\,28\}$ is congruent to one of the representative days of the month $1$ to $7$. The $7\ssptimes 7$ sub-matrix is a symmetric {\sl Toeplitz} matrix with maximum $688$ only in the main diagonal. For $d\sspeq 29$ the maximum $644$ appears for the days of the week $D\sspeq 0$ and $2$. For $d\sspeq 30$ the maximum $631$ appears for $D\sspeq 1$ and $3$, and for $d\sspeq 31$ the maximum $402$ appears for $D\sspeq 4$. 
\section{Multiplicities for days of the month on days of the week}
\hskip 1cm 
The focus of this note is on the number of times a given day of the month                         
falls on a day of the week for each of the $400$ years of a {\sl Gregorian}                        
cycle. That is we wish to know, for example, how many Friday $13$th                         
happen in a given year, not only the total sum of such events over a $400$                         
year cycle. To this end we consider sequences $M(d,\,D)$ of length $400$                     for the consecutive years $y\sspeq  0,\,1,\,...,\,399$, with entry at position $y$ giving the number of times an admissible day of the month $d\sspin \{1,\,2,\,..,\,31\}$ falls on a day of the week $D \in\{0,1,...,6\}$. (This $M$ should not be confused wit the $M$ of eq.~\ref{(1.1)}. $M(d,\, D)$ is then repeated periodically for values $y\sspgeq 400$. \pn
The representatives of the conjugacy classes modulo $7$ for the days $d \sspin \{1,\,2,\,...,\,28\}$ are  $d \sspin \{1,\,2,\,...,\,7\}$. $d\sspeq 1$ represents the days from $\{1,\, 8,\, 15,\, 22\}$, {\it etc.} For every $d\sspin \{2,\,...,\,28\}$ the sequence $M(d,\,D)$ can be obtained from one of the seven  $M(1,\,D)$ sequences. The precise relation is $M(d,\,D) \sspeq M(\bar d,\,D) \sspeq M(1, (D\sspp 1 \sspm \bar d)$ mod $7)$, with $\bar d\sspdef 7$ if $d$\,mod $7\sspeq 0$ and otherwise $d$\,mod $7$. This is exhibited as a  $7\times 7$ {\sl Toeplitz} type matrix scheme in {\it Table} $2$.
\pn
Therefore, only $4\cdot 7\sspeq 28$ such sequences of length $400$ are needed, viz $M(1,\,D), M(29,\,D), M(30,\, D)$ and $M(31,\, D)$, for $D\sspin {0,\,1,\,...,\,6}$. 
\pn
The main observation reported here is that each of these $28$ basic sequences                           
$M(d,\,D)$ of lenght $400$, used periodically, can be put together from pieces of one of only $4$ cyclic sequences of length $28$. These four short sequences will be called $S_{1}$, $S_{29}$, $S_{30}$ and $S_{31}$. They will be used periodically to build $M(1,\,D), M(29,\,D), M(30,\, D)$ and $M(31,\, D)$, respectively. They are given in the second row of {\it Table} $3$. Note that the offset is $0$ for these sequences.
\begin{eqnarray}
S_1\sspeq    &[1,\,2,\,2,\,1,\,2,\,1,\,2,\,2,\,1,\,3,1,\,1,\,3,\,2,\,1,3,\,1,\,2,\,2,\,2,\,2,\,1,\,1,\,2,\,2,\,1,\,3,\,1], \label{(2.1)}&\\ 
S_{29}\sspeq &[1,\,2,\,2,\,1,\,2,\,1,\,2,\,2,\,1,\,2,\,1,\,1,\,3,\,2,\,1,\,2,\,1,\,2,\,2,\,2,\,2,\,1,\,1,\,2,\,2,\,1,\,2,\,1],\label{(2.2)}&\\
S_{30}\sspeq &[3,\,2,\,1,\,2,\,1,\,2,\,2,\,2,\,2,\,1,\,1,\,2,\,2,\,1,\,2,\,1,\,1,\,2,\,2,\,1,\,1,\,1,\,2,\,2,\,1,\,2,\,1,\,1],\label{(2.3)}&\\
S_{31}\sspeq &[1,\,0,\,1,\,1,\,1,\,1,\,1,\,0,\,1,\,1,\,2,\,1,\,0,\,1,\,1,\,1,\,2,\,1,\,0,\,1,\,1,\,2,\,1,\,1,\,1,\,1,\,1,\,2].\label{(2.4)}&
\end{eqnarray}
These four sequences for multiplicities are obtained from the first $28$ entries of $M(1,\,0),\,M(29,\,0),\,M(30,\, 0)$ and $M(31,\, 0)$, respectively (also with offset $0$).\pn
These sequences are used cyclically and shown as 'clocks' in the {\sl Figure}s $1,\,2,\,3$ and $4$.\pn
Later they will be considered without the outer brackets and denoted by $\bar S $.
\pn  
The composition of the $28$ sequences $M(i,\,D)$ for $i\sspeq 1,\, 29,\, 30,\, 31$ and $D\sspeq 0,\,1,\, ...\,6$, are encoded like shown in {\it Table} $3$. This encoding is explained by giving examples for a certain day of the week $D$ for each $M(i,\, D)$ in the following four subsections.
\pbn
{\meinfont{\bf A) Example $\bf M(1,\,2)$}}
\psn
Here the cyclic sequence with period $S_1$ is used. For $D\sspeq 2$ (Tuesday) the start in the representative year $y\sspeq 0$ is seen from the number of entries $2$ of the list $[D(1,m,0)_{m=1..12}]\sspeq [6,\, 2,\, 3,\, 6,\, 1,\, 4,\, 6,\, 2,\, 5,\, 0,\, 3, \,5]$ to be $M(1,\,2)(0) \sspeq 2$. For the next years the multiplicity of 2 is found by starting with $S_1(4)\sspeq 2$ for the all together $100\sspeq 3\cdot28 \sspp 16$ entries, repeating $S_1$, shown in {\sl Figure} $1$ cyclically. This means on starts with $S_1(4)$ reaching first the end of the cycle $S_1(27)$, then repeating the whole $S_1$ twice, and in order to obtain $100 \sspeq (27\sspm 4\sspp 1) \sspp 2\cdot 28 \sspp (4 \sspp 16) $ entries, one keeps going for the missing $4 \sspp 16\sspeq 20$ entries to end up with $S_1(19)\sspeq 2$. Here one has to stop because the next entries in $M(1,\,2)$ are $1,\,3,\,1,\,1,\,3 \,...$, but $S_1(20)\sspeq 2$, not $1$.\pn
This second piece of $M(1,\,2)$ starts with $S_1(8)$, \ie one has to move from $S_1(19)$ to $S_1(8)\sspeq 1$ in $17$ steps. Up to now the encoding is $(4)100(17)$, where the numbers in brackets give the number of steps one has to move from the reached entry of $S_1$ to the start of the next valid sub-sequence read off the $S_1$ clock ($(4)$ indicates the four steps from $S_1(0)$ to $S_1(4)$), and the following numbers give the length of the sequence cycling in $S_1$, before the next wrong number compared to $M(1,\,2)$ appears. This second sub-sequence has $101 \sspeq (27-8+1) \sspp 2\cdot 28\sspp (8+17)$ entries, viz $S_1(8)$ to  $S_1(27)$, then cycling $S_1$ twice and going fom $S_1(0)$ to $S_1(24)\sspeq 2$. The next entry $S(25)\sspeq 1$ does not fit $M(1,\,2)(201)\sspeq 2$, which turns out to be $S_1(13)$, after again $17$ steps.
Up to now the encoding is $(4)100(17)101(17)$.
\pn
The third piece of $M(1,\,2)$ has $102\sspeq (27-13+1)\sspp 3\cdot 28\sspp 3$ entries, viz $S_1(13)$ to $S_1(27)$, then cycling thrice, and continuing with  $S_1(0)$ to $S_1(2)$. The next entry $S_1(3)\sspeq 1$ does not fit $M(1,\,2)(303)\sspeq 2$. The correct fourth sub-sequence starts with $S_1(19)\sspeq 2$, after again $(17)$ steps. The encoding continues therefore with $(4)100(17)101(17)102(17)$.
\pn
The last piece of  $M(1,\,2)$ has the missing $97\sspeq (27-19+1)\sspp 3\cdot 28\sspp 4$ entries, viz $S_1(19)$ to $S_1(27)$, cycling thrice and going from $S_1(0)$ to $S_1(3)$\sspeq 1. Then a new length $400$ cycle will starts for $M(1,\, 2)$ with $S_1(4) \sspeq 2$ (after one step).\
Thus the complete encoding of this length $400$ sequence is
\Beq
{\rm encoded}\ M(1,\,2) \sspeq (4)100(17)101(17)102(17)97,     \label{(2.5)}
\Eeq
as given in {\it Table} $3$.
\pn
Written in terms of the four $S_1$ pieces of lengths $100, 101, 102$ and $97$ this becomes 
\Beqarray
M(1,\,2)\sspeq && [S_1(4..27),\bar S_1,\bar S_1,S_1(0..19),\ \ S_1(8..27),\bar S_1,\bar S_1,S_1(0..24),\ \ S_1(13..27)\bar S_1,\bar S_1,\bar S_1,S_1(0..2), \nonumber \\    
&& S_1(19..27)\bar S_1,\bar S_1,\bar S_1,S_1(0..3)], \label{(2.6)}
\Eeqarray
where $S_1(i..j)$ stands for the sub-sequence $S_1(k)_{k=i,i+1,...,j}$ (with entries separated by commas), and $\bar S_1$ is the sequence $S_1$ given in eq.~\ref{(2.1)} without the brackets [.].
\pbn
\vfill\eject\noindent
{\meinfont{\bf B) Example $\bf M(29,\,1)$}}
\psn
After the detailed description of the previous example it suffices to give the encoding, from {\it Table} $3$ and the corresponding decoding in terms of five pieces of the sequence $S_{29}$, shown in {\it Fig.} $2$.
\Beq
\rm{encoded}\,M(29,\,1) \sspeq (16)102(17)98(17)100(7)4(11)96\,. \label{(2.7)}
\Eeq
\Beqarray
M(29,\,1) \sspeq && [S_{29}(16..27),\bar S_{29},\bar S_{29},\bar S_{29},S_{29}(0..5),\ \ S_{29}(22..27),\bar S_{29},\bar S_{29}\bar S_{29},S_{29}(0..7),\ \ \nonumber \\
&&  S_{29}(24..27),\bar S_{29},\bar S_{29},\bar S_{29},S_{29}(0..11),\ \  S_{29}(18..21),\ \ S_{29}(4..27)\bar S_{29},\bar S_{29},S_{29}(0..15)]\,. \label{(2.8)}
\Eeqarray
\pbn
{\meinfont{\bf C) Example $\bf M(30,\,4)$}}
\psn
This case consists of five pieces, see {\it Table} $3$, and $S_{30}$ is shown in {\it Fig.} $3$.
\psn
\Beq
\rm{encoded}\,M(30,\,4) \sspeq (8)100(17)100(7)5(11)96(17)99\,. \label{(2.9)}
\Eeq
\Beqarray
M(30,\,4) \sspeq && [S_{30}(8..27),\bar S_{30},\bar S_{30}, S_{30}(0..23),\ \ S_{30}(12..27),\bar S_{30},\bar S_{30},\bar S_{30},\ \ S_{30}(6..10),\ \ \nonumber \\
&& S_{30}(21..27),\bar S_{30},\bar S_{30},\bar S_{30},S_{30}(0..4),\ \  S_{30}(21..27),\bar S_{30},\bar S_{30},\bar S_{30},S_{30}(0..7)]. \label{(2.10)}
\Eeqarray
\pbn
{\meinfont{\bf D) Example $\bf M(31,\,0)$}}
\psn
This case consists of four pieces, see {\it Table} $3$, and $S_{31}$ is shown in {\it Fig.} $3$.
\psn
\Beq
\rm{encoded}\,M(31,\,0) \sspeq (0)102(17)99(17)99(17)(100)\,. \label{(2.11)}
\Eeq
\Beqarray
M(31,\,0) \sspeq && [\bar S_{31},\bar S_{31},\bar S_{31},S_{31}(0..17),\ \ S_{31}(6..27),\bar S_{31},\bar S_{31},S_{31}(0..20),\ \ \nonumber \\
&& S_{31}(9..27),\bar S_{31},\bar S_{31},S_{31}(0..23),\ \ S_{31}(12..27),\bar S_{31},\bar S_{31},\bar S_{31}]. \label{(2.12)}
\Eeqarray
\pbn                                                                                         
                                                                           
\pbn
\pbn
\hrule
\psn
{\bf Keywords:}: Gregorian calendar, days of month on days of the week, integer sequences of period $28$. \psn
{\bf MSC-numbers:} 11N69, 11Y55, 05Axx
\psn
\hrule
\vfill                                                                                            
\eject
\pbn
\begin{center}
{\includegraphics[height=15cm,width=.8\linewidth]{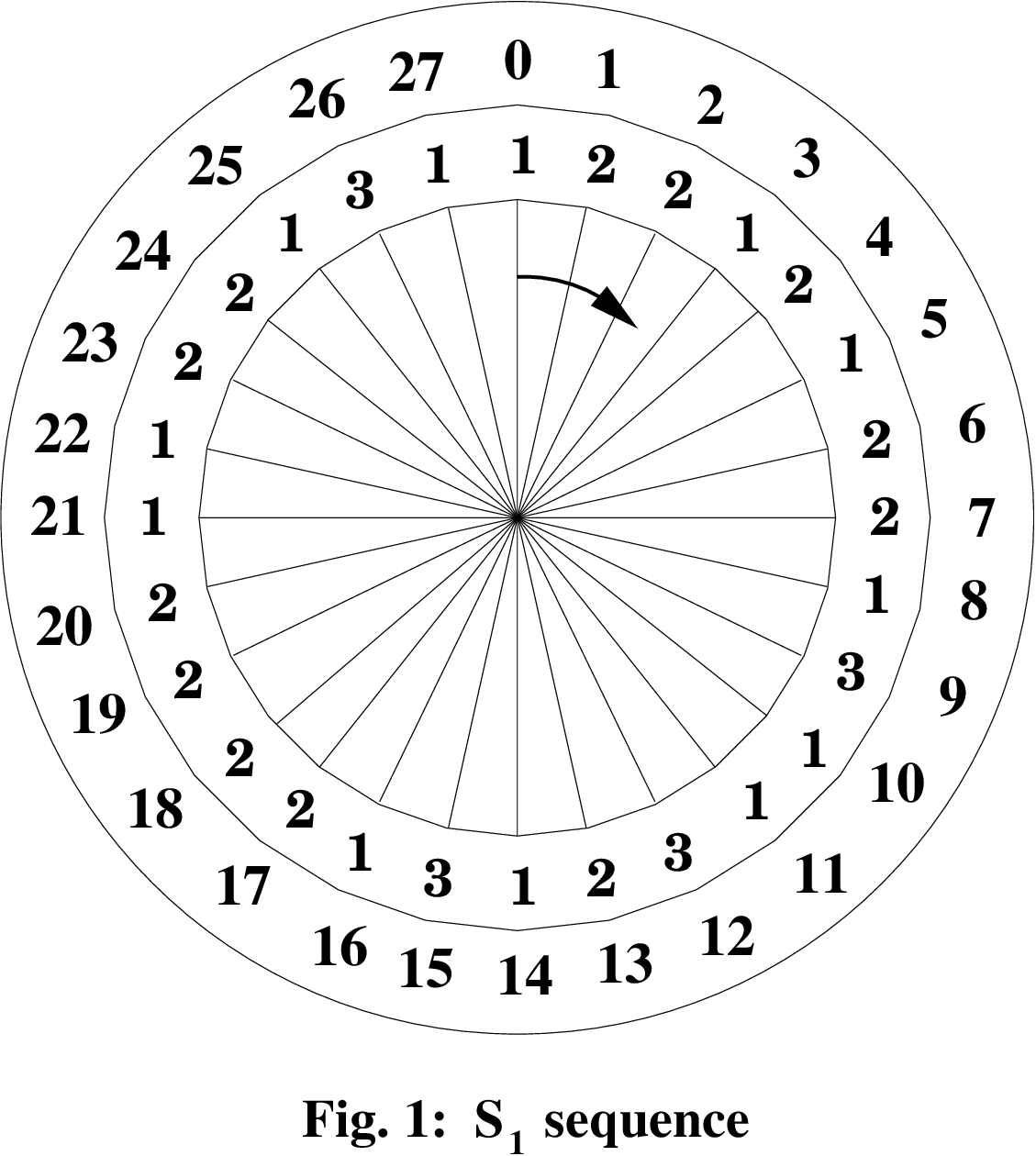}
{\hfill}}
\end{center}
\pbn
\vfill\eject
\begin{center}
{\includegraphics[height=15cm,width=.8\linewidth]{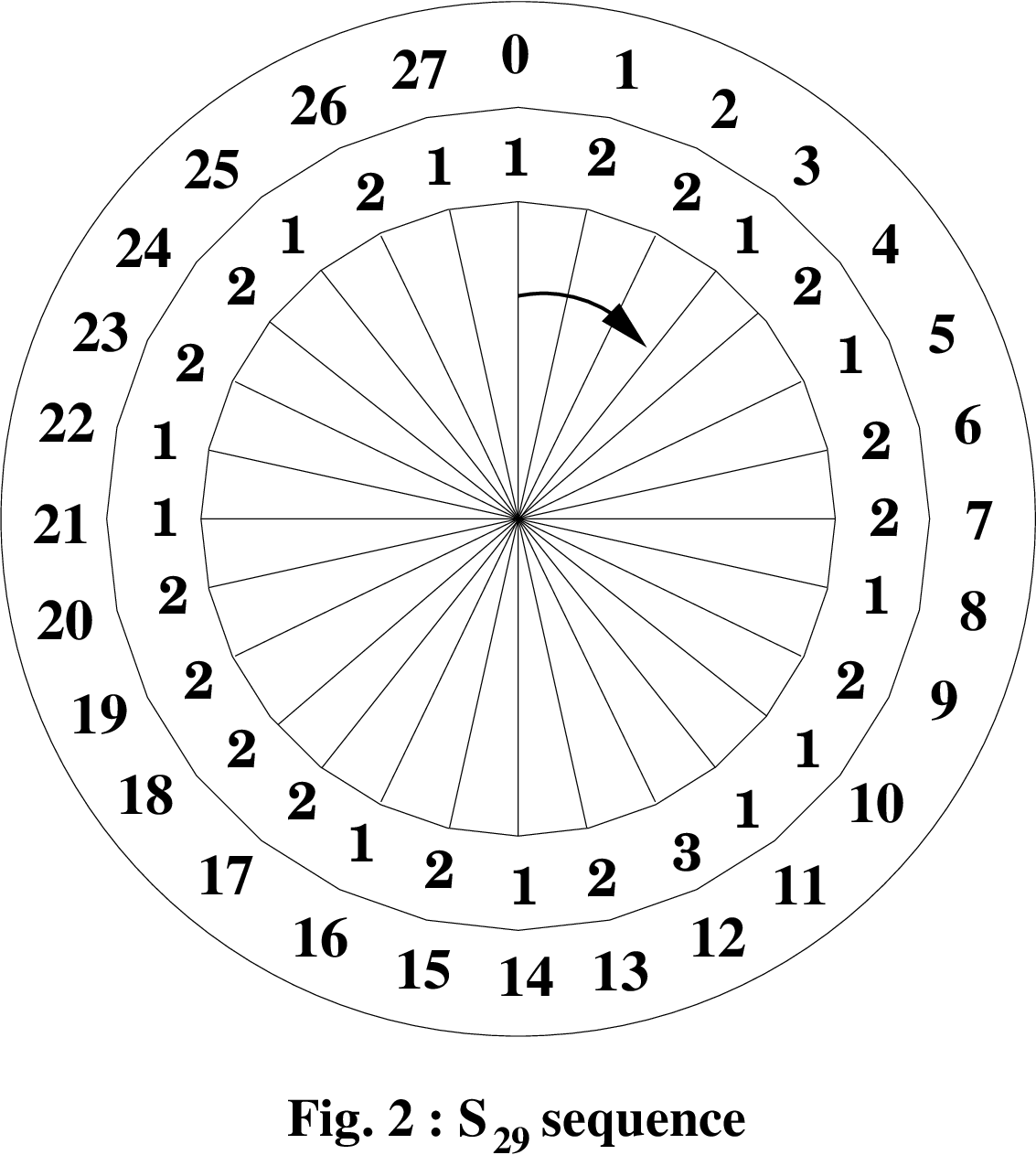}
{\hfill}}
\end{center}
\pbn
\vfill\eject
\begin{center}
{\includegraphics[height=15cm,width=.8\linewidth]{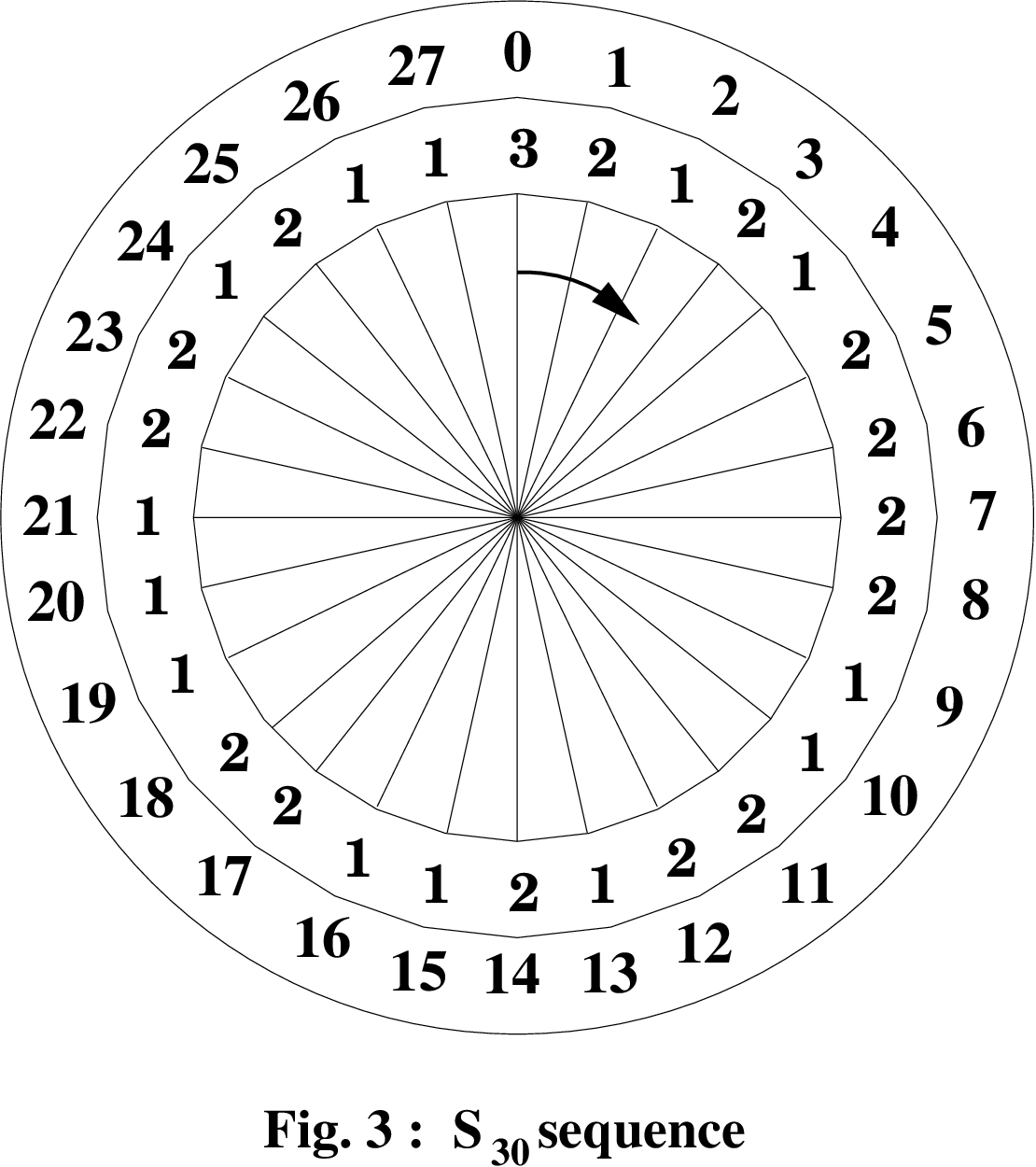}
{\hfill}}
\end{center}
\psn
\vfill\eject
\begin{center}
{\includegraphics[height=15cm,width=.8\linewidth]{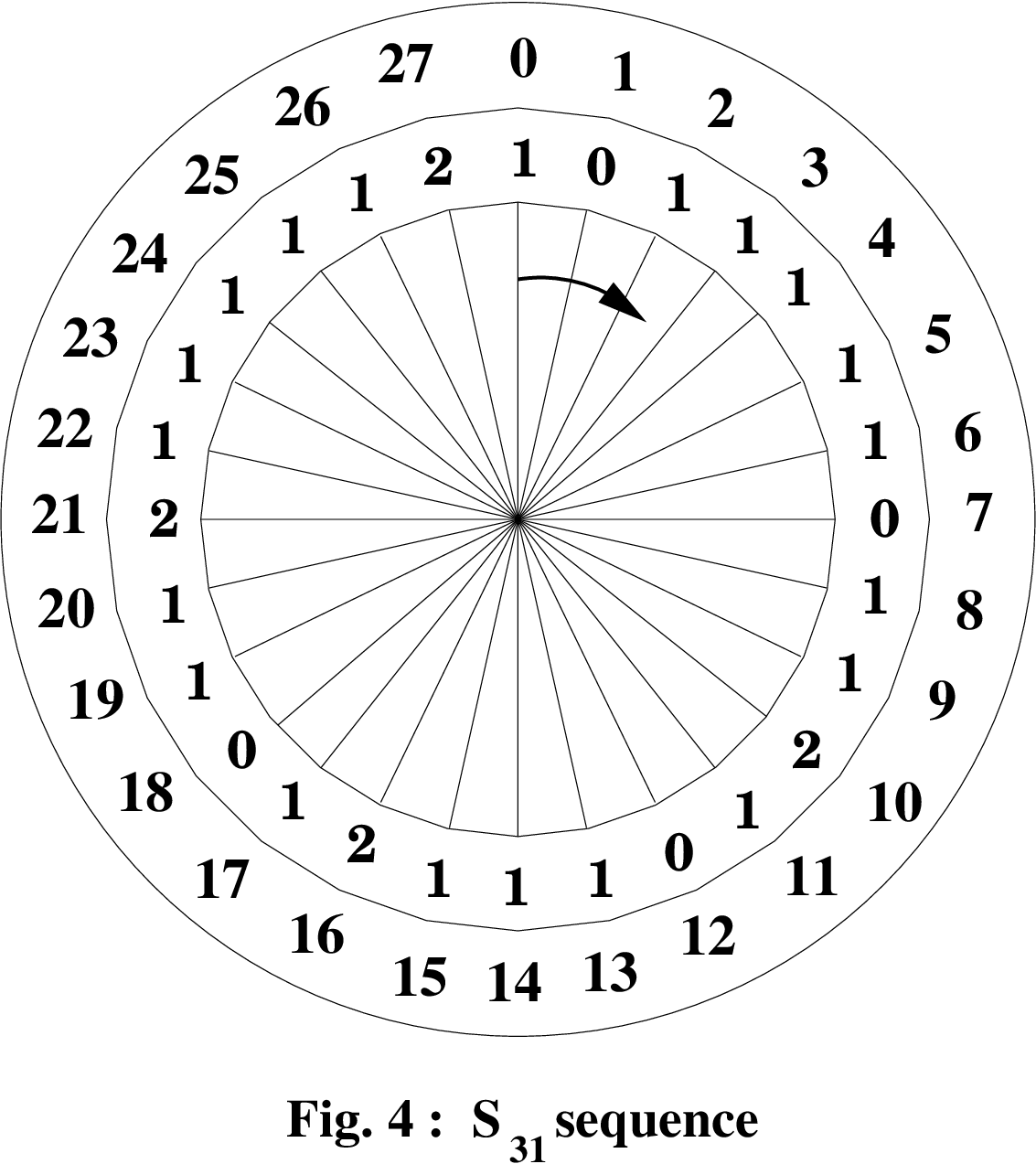}
{\hfill}}
\end{center}
\pbn
\vfill \eject
\begin{center}
{\large {\bf Table 1: Occurrences for days of a month $\bf d$ at a weekday $\bf D$\\ 
for each Gregorian cycle of $\bf 400$ years\phantom{xxxx}}}
\par\bigskip\noindent
\begin{tabular}{|l||c|c|c|c|c|c|c|}\hline
&& && && & \\
$\bf{d \backslash D}$ &  $\bf 0$ & $\bf 1$ & $\bf 2$ & $\bf 3$  &  $\bf 4$ & $\bf 5$ & $\bf 6$ \\
   &  $\bf Su$ & $\bf Mo$ & $\bf Tu$ & $\bf We$  &  $\bf Th$ & $\bf Fr$ & $\bf Sa$ \\
&& && && & \\ 
\hline\hline 
$\bf 1\,(\ \, 8,15,22)$ &  {\color{magenta}$\bf 688 $}  & $\bf 684 $ & $\bf 687 $  &  $\bf 685 $& $ \bf 685 $ & $ \bf 687 $ & $\bf 684$\\
$\bf 2\,(\ \,9,16,27)$ &  $\bf 684 $  & {\color{magenta} $\bf 688$} & $\bf 684 $  &  $\bf 687 $  & $\bf 685 $ & $\bf 685 $ & $\bf 687 $\\
$\bf 3\,(10,17,24)$ &  $\bf 687 $  & $\bf 684$ & {\color{magenta}$\bf 688 $}  &  $\bf 684  $ & $\bf 687 $ &  $\bf 685 $ & $\bf 685 $\\
$\bf 4\,(11,18,25)$ &  $\bf 685 $  & $\bf 687 $ & $\bf 684 $  &  {\color{magenta}$\bf 688 $} & $\bf 684 $ &  $\bf 687 $ & $\bf 685 $\\
$\bf 5\,(12,19,26)$ &  $\bf 685 $  & $\bf 685 $ & $\bf 687$  &  $\bf 684 $ & {\color{magenta} $\bf 688 $} &  $\bf 684 $ & $\bf 687$\\
$\bf 6\, (13,20,27)$ &  $\bf 687 $  & $\bf 685$ & $\bf 685$  &  $\bf 687 $& $\bf 684 $ &  {\color{magenta} $\bf 688 $} & $\bf 684$\\
$\bf 7\,(14,21,28)$ &  $\bf 684 $  & $\bf 687$ & $\bf 685$  &  $\bf 685 $& $\bf 687 $ &  $\bf 684 $ &{\color{magenta} $\bf 688$}\\
$\bf 29$ & {\color{magenta} $\bf 644 $}  & $\bf 641$ & {\color{magenta}$\bf 644$}  &  $\bf 642 $& $\bf 642 $ &  $\bf 643 $ & $\bf 641$\\
$\bf 30$ & $\bf 627 $  & {\color{magenta}$\bf 631 $} & $\bf 626 $  & {\color{magenta} $\bf 631 $}& $\bf 627 $ &  $\bf 629 $ & $\bf 629 $\\
$\bf 31$ & $\bf 400 $  & $\bf 399 $ & $\bf 401 $  & $\bf 398 $ & {\color{magenta}$\bf 402 $} &  $\bf 399 $ & $\bf 401$\\
\hline
\hline
\end{tabular}
\end{center}
\pbn
\begin{center}
{\large {\bf Table 2: \  $\bf T_7(\overline d,\,D)$ matrix for the identity 
\phantom{xxxxxx}\\ \hskip 3.4cm $\bf M(d,\,D) \sspeq  M(1,\,T_7(\overline d,\ D))$ for $\bf d\sspin \{1,\,2,\,...,28\}$}}
\par\bigskip\noindent
\begin{tabular}{|l||c|c|c|c|c|c|c|}\hline
&& && && & \\
$\bf{\overline d \backslash D}$ &  $\bf 0$ & $\bf 1$ & $\bf 2$ & $\bf 3$  &  $\bf 4$ & $\bf 5$ & $\bf 6$ \\
&& && && & \\ 
\hline\hline 
$\bf 1$ &  $\bf 0 $ & $\bf 1 $ & $\bf 2 $  &  $\bf 3 $& $ \bf 4 $ & $ \bf 5 $ & $\bf 6$\\
$\bf 2$ &  $\bf 6 $ & $\bf 0 $ & $\bf 1 $  &  $\bf 2 $& $ \bf 3 $ & $ \bf 4 $ & $\bf 5$\\
$\bf 3$ &  $\bf 5 $ & $\bf 6 $ & $\bf 0 $  &  $\bf 1 $& $ \bf 2 $ & $ \bf 3 $ & $\bf 4$\\
$\bf 4$ &  $\bf 4 $ & $\bf 5 $ & $\bf 6 $  &  $\bf 0 $& $ \bf 1 $ & $ \bf 2 $ & $\bf 3$\\
$\bf 5$ &  $\bf 3 $ & $\bf 4 $ & $\bf 5 $  &  $\bf 6 $& $ \bf 0 $ & $ \bf 1 $ & $\bf 2$\\
$\bf 6$ &  $\bf 2 $ & $\bf 3 $ & $\bf 4 $  &  $\bf 5 $& $ \bf 6 $ & $ \bf 0 $ & $\bf 1$\\
$\bf 7$ &  $\bf 1 $ & $\bf 2 $ & $\bf 3 $  &  $\bf 4 $& $ \bf 5 $ & $ \bf 6 $ & $\bf 0$\\
\hline
\hline
\end{tabular}
\end{center}
\psn
\begin{center}
\bf Example: $ \bf M(17,\,4)\sspeq M(\overline {17},\, 4)\sspeq M(3,\,4)\sspeq M(1,\, 2)$ \\
\hskip 3.9cm with $\bf \overline d\sspdef 7$\ if $\bf d\sspequiv 0$\,(mod$\, 7)$, otherwise $\bf d$\,mod \,$7$. 
\end{center}
\pbn
\pbn
\pbn
\begin{landscape}
\begin{center}
{\large {\bf Table 3: \ Sequences for number of times a given days of the month $\bf d$ falls on days of a week within a $\bf 400$ cycle}}
\par\bigskip\noindent
\begin{tabular}{|l||l|l|l|l|}\hline
&& && \\ &  \hskip 1.5cm $\bf d \sspin\{1,\,8,\, 15,\, 22\}$ & \hskip 2.5cm$\bf d\sspeq 29 $ & \hskip 2.5cm$\bf d\sspeq 30$ & \hskip 2.5cm$\bf d\sspeq 31$\\
&& &&  \\ 
\hline\hline
&&&& \\
\bf Period 28 & $\bf S_1\sspeq [1,2,2,1,2,1,2,2,1,3,1,1, $ & $\bf S_{29}\sspeq [1,2,2,1,2,1,2,2,1,2,1,1,$ & $ \bf S_{30}\sspeq [3,2,1,2,1,2,2,2,2,1,1,2,$ & $\bf S_{31}\sspeq [1,0,1,1,1,1,1,0,1,1,2,1,$\\
\bf sequence & $ \bf 3,2,1,3,1,2,2,2,2,1,1,2,2,1,3,1]$ & $ \bf 3,2,1,2,1,2,2,2,2,1,1,2,2,1,2,1]$  & $ \bf 2,1,2,1,1,2,2,1,1,1,2,2,1,2,1,1$& $\bf 0,1,1,1,2,1,0,1,1,2,1,1,1,1,1,2]$\\
&&&& \\
\hline
\hline
\bf Sunday $(\bf 0)$ & $(0)101(17)99(17)101(17)99 $ & $(0)101(17)99(17)101(17)99 $ &  $ (0)103(17)97(17)100(17)100 $ & $ (0)102(17)99(17)99(17)100$\\
\hline
\bf Monday $(\bf 1)$ & $(16)100(23)6(23)94(17)100(7)4(11)96 $ & $ (16)102(17)98(17)100(7)4(11)96 $ & $(16)101(17)99(17)101(17)99 $ & $ (16)101(17)101(17)98(17)100$\\
\hline
\bf Tuesday $(\bf 2)$ & $(4)100(17)101(17)102(17)97$ & $(4)100(17)101(17)102(17)97 $ & $ (4)100(17)100(17)100(7)5(11)95 $ & $(4)101(17)99(17)101(17)99$ \\
\hline
\bf Wednesday $(\bf 3)$ & $(20)100(17)100(7)4(11)97(17)99 $ & $(20)100(17)100(7)4(11)97(17)99 $ & $(20)100(17)101(17)102(17)97 $ & $(20)102(17)98(17)102(17)98 $ \\
\hline
\bf Thursday $(\bf 4)$ & $ (8)101(17)102(17)97(23)6(23)94 $ & $(8)101(17)102(17)99(17)98 $ & $(8)100(17)100(7)5(11)96(17)99$ & $(8)100(17)101(17)100(17)99 $ \\
\hline
\bf Friday $(\bf 5)$ & $(24)100(7)4(11)97(17)99(17)100 $ & $(24)100(7)4(11)97(17)100(16)99 $ & $(24)101(17)102(17)97(17)100$ & $(24)100(17)102(17)99(17)99 $ \\
\hline
\bf Saturday $(\bf 6)$ & $(12)103(17)97(23)6(23)94(17)100 $ & $(12)103(17)99(17)103(18)95 $ & $(12)100(7)5(11)96(17)100(16)99 $ & $(12)101(17)100(17)101(17)98$ \\
\hline
\end{tabular}
\end{center}
\end{landscape}
\pbn
\end{document}